\documentclass{amsart}
\usepackage{latexsym}
\usepackage{amsmath}
\usepackage{amsfonts}
\usepackage{amssymb}
\usepackage{amsthm}

\theoremstyle{plain}
\newtheorem{theorem}{Theorem}[section]
\newtheorem{proposition}[theorem]{Proposition}
\newtheorem{lemma}[theorem]{Lemma}
\newtheorem{corollary}[theorem]{Corollary}

\theoremstyle{definition}
\newtheorem{example}[theorem]{Example}
\newtheorem{definition}[theorem]{Definition}
\newtheorem{remark}[theorem]{Remark}

\newcommand\bt{\begin{theorem}}
\newcommand\et{\end{theorem}}

\newcommand\bl{\begin{lemma}}
\newcommand\el{\end{lemma}}

\newcommand\bp{\begin{proposition}}
\newcommand\ep{\end{proposition}}

\newcommand\bc{\begin{corollary}}
\newcommand\ec{\end{corollary}}

\newcommand\bd{\begin{definition}}
\newcommand\ed{\end{definition}}

\newcommand\br{\begin{remark}}
\newcommand\er{\end{remark}}

\newcommand\bex{\begin{example}}
\newcommand\eex{\end{example}}

\newcommand\bpf{{\it Proof. }}
\newcommand\epf{$\,$ $\Box$ \\}

\title[ON FINITE MOLECULARIZATION DOMAINS]{ON FINITE MOLECULARIZATION DOMAINS}

\author[Hetzel]{Andrew J. Hetzel}
\address{Department of Mathematics, Tennessee
Tech University, 110 University Drive, Bruner Hall 235, Cookeville, Tennessee 38505, USA}
\email{ahetzel@tntech.edu}

\author[Lawson]{Anna L. Lawson}
\address{Department of Mathematics, The University of Tennessee, 227 Ayres Hall, 1403 Circle Drive, Knoxville, TN 37996, USA}
\email{alitchfo@vols.utk.edu}

\author[Reinhart]{Andreas Reinhart}
\address{Department of Mathematical Sciences, New Mexico State University, 1290 Frenger Mall, Las Cruces, New Mexico 88003, USA}
\email{andreas.reinhart@uni-graz.at}

\thanks{This work is based in part on the second-named author's
master's research at Tennessee Tech University.}

\thanks{The third-named author of this work was supported by the Austrian Science Fund FWF, Project Number J4023-N35.}

\thanks{The authors wish to express their thanks to the referee of this article for his thorough and thoughtful reading of the manuscript.}


\begin{document}

\begin{abstract}
In this paper, we advance an ideal-theoretic analogue of a ``finite factorization domain" (FFD), giving such a domain the moniker ``finite molecularization domain" (FMD). We characterize FMD's as those factorable domains (termed ``molecular domains" in the paper) for which every nonzero ideal is divisible by only finitely many nonfactorable ideals (termed ``molecules" in the paper) and the monoid of nonzero ideals of the domain is unit-cancellative, in the language of Fan, Geroldinger, Kainrath, and Tringali. We develop a number of connections, particularly at the local level, amongst the concepts of ``FMD", ``FFD", and the ``finite superideal domains" (FSD's) of Hetzel and Lawson. Characterizations of when $k[X^2, X^3]$, where $k$ is a field, and the classical $D+M$ construction are FMD's are provided. We also demonstrate that if $R$ is a Dedekind domain with the finite norm property, then $R[X]$ is an FMD.
\end{abstract}

\maketitle

\section{INTRODUCTION} \label{S1}

Throughout this paper, all rings are commutative with $1 \neq 0$. For over a century, the study of various types of decompositions of ideals has occupied an important place in commutative ring theory. While the famed Noether-Lasker theorem is often considered the archetype for such investigations, invaluable research has been conducted related to decomposing ideals as a \textit{product} of a certain type of ideal; see, amongst many others, \cite{MS},\cite{Olb}, \cite{VY},\cite{AS},\cite{AM} (it should be noted that, in general, the term ``factoring" is used with regard to products of ideals while the term ``decomposing" is used with regard to intersections of ideals). Moreover, a modern focus on products of ideals has tremendous worth even in the classical context of algebraic geometry that motivated the Noether-Lasker theorem. To wit, in light of the fact that the algebraic variety of a product of ideals is the same as the variety of the corresponding intersection, computing a basis for a product of ideals is far more straightforward than computing a basis for an intersection of ideals.

In addition, while the idea of an ``irreducible ideal" is quite natural for considerations of decomposing ideals, it was not until 1964 that a truly parallel notion, embodied by the concept of a ``nonfactorable ideal", for factoring ideals was advanced by H.S. Butts \cite{HSB}. A {\em nonfactorable ideal} $I$ of a commutative ring $R$ is a nonzero proper ideal of $R$ such that whenever $I = JK$ for some ideals $J$ and $K$ of $R$, it must be the case that either $J = R$ or $K = R$ (see also \cite{HSB2}). In \cite{HSB}, Butts demonstrated that if $R$ is an integral domain, then every nonzero, proper ideal of $R$ can be factored uniquely (up to the order of the factors) as a product of nonfactorable ideals of $R$ if and only if $R$ is a Dedekind domain. As such, the concept of ``Dedekind domain" can be viewed as the proper ideal-theoretic analogue of ``unique factorization domain".

Arguably inspired by the tremendous fruitfulness in studying certain generalizations of a ``unique factorization domain" in \cite{AAZ1}, D.F. Anderson, H. Kim, and J. Park \cite{AKP} introduced and explored {\em factorable domains}--domains $R$ with the property that every nonzero proper ideal of $R$ is a product of nonfactorable ideals of $R$--an ideal-theoretic analogue to atomic domains. Consistent with this perspective, in this paper, we wish to advance an ideal-theoretic analogue of the notion of a ``finite factorization domain" (or FFD), one of the generalizations of a ``unique factorization domain" introduced in \cite{AAZ1}. On our way to discovering definitive information for such an analogue in the contexts of certain polynomial rings (see Theorem \ref{theo10}, Corollary \ref{coro14.1}) and the classical $D+M$ construction (see Theorem \ref{theo16}), we pick up some novel information about factorable domains (notably Theorem \ref{theo15}) and even nonfactorable ideals themselves (notably Proposition \ref{prop2}, Proposition \ref{prop11}, Proposition \ref{prop12}, Theorem \ref{theo13}).

At this point, a major caveat concerning terminology in this paper is warranted. In spite of the use in \cite{HSB}, \cite{HSB2}, \cite{AKP}, \cite{HL1}, and \cite{L} of the terms ``nonfactorable ideal", ``factorable domain", and ``factorable ring", we do not wish to continue this practice in the present paper. The main reasons for this are that (1) the latter two terms create an ambiguity if taken outside of context, as there are a fair number of different types of factorizations of ideals (as mentioned above) that a term like ``factorable" may reference, and (2) such terms do not seem to connect with the element-level inspiration for these notions, where expressions such as ``atom" and ``atomic domain" are standard. To begin remedying these issues, we have chosen to adopt the terms ``molecule" for ``nonfactorable ideal", ``molecular domain" for ``factorable domain", and ``molecularization" to mean ``product of molecules". Moreover, the term ``finite ideal factorization domain" in \cite{L} will henceforth be replaced with ``finite molecularization domain", the focal concept of this paper. In addition, a nonzero proper ideal that is not a molecule will be called ``compound". Overall, such a change in vocabulary has the virtue of creating a certain idiomatic aesthetic, particularly in view of results such as Proposition \ref{prop2}, Corollary \ref{coro2.5}, and Proposition \ref{prop9} in this paper.

As usual, the set of whole numbers (that is, the set of nonnegative integers) will be represented by $\mathbb{W}$. Let $R$ be a domain. The group of units of $R$ will be designated by $U(R)$. Throughout this paper, the {\em dimension} of $R$, denoted $\mbox{{dim}}(R)$, always refers to the Krull dimension of $R$, that is, the supremum of the lengths of all chains of prime ideals of $R$. In particular, if $P$ is a prime ideal of $R$, then the {\em height} of $P$ is the dimension of the localization $R_P$. If $I$ is an ideal of $R$, then the nilradical of $I$ is given by $\sqrt{I} = \cap \{P \mid P \textrm{ is a prime ideal of } R \textrm{ containing } I\}$. If $I$ and $J$ are ideals of $R$, we say that $J$ {\em divides} $I$ if there exists an ideal $K$ of $R$ for which $I = JK$.

Furthermore, we distinguish between calling the domain $R$ {\em quasilocal} if it has a unique maximal ideal and {\em local} if additionally $R$ is Noetherian. Similarly, we distinguish between calling $R$ {\em quasisemilocal} if it has only finitely many maximal ideals and {\em semilocal} if additionally $R$ is Noetherian. The integral closure of $R$ in some (given) field extension of the quotient field of $R$ will be denoted by $\overline{R}$; where specification is not provided, the field extension may be assumed to be the quotient field itself. An {\em overring} of $R$ is a ring containing $R$ and contained within the quotient field of $R$. As in \cite{HL0}, $R$ is called a {\em finite superideal domain} (or FSD) if every nonzero proper ideal of $R$ has only finitely many ideals of $R$ containing it.

Any unexplained terminology is standard, as in \cite{Gil}, \cite{Kap}, \cite{AtM}.

\section{PROPERTIES OF MOLECULES} \label{S2}

We begin this paper by presenting several novel properties of molecules. Our inaugural results, Proposition \ref{prop2} below and its associated Corollary \ref{coro2.5}, draw valuable connections between the type of ideal-level factorizations being explored in this paper and the corresponding type of element-level factorizations as a product of irreducibles. In particular, the relationship between molecules and atoms is most intimate for those domains with trivial Picard groups (which include the respective classes of quasisemilocal domains, B\'{e}zout domains, UFD's, and, thanks to \cite[Theorem 6.1]{Olb2}, one-dimensional domains with nonzero Jacobson radical).

\bp \label{prop2}
Let $R$ be a domain and $I=(a)$ a principal ideal of $R$. If $I$ is a molecule, then $a$ is an atom. If $R$ is further assumed to have trivial Picard group, then the converse is true, as well. Moreover, in this context, if $J$ is an ideal that divides the principal ideal $I$, then $J$ is also principal.
\ep

\bpf
Let $R$ be a domain and let $I=(a)$ be a principal ideal of $R$. If $a$ is not an atom, then $a = bc$ for some nonunits $b$ and $c$, which implies $I = (b)(c)$, where $(b)$ and $(c)$ are proper ideals of $R$. Hence, $I$ is not a molecule.

Now, further suppose that $R$ has trivial Picard group and that $a$ is an atom. Assume to the contrary that $I = JK$ for proper ideals $J$ and $K$. Since $J$ and $K$ are necessarily invertible, it must be the case that $J$ and $K$ are principal. Thus, $I = (b)(c)$, and so $a = ubc$ for some unit $u$. However, the assumption that $J$ and $K$ are proper guarantees that $b$ and $c$ are nonunits, contradicting the irreducibility of $a$. Therefore, it must be the case that $I$ is a molecule.

Finally, observe that even without the assumption that $a$ is an atom, the above work shows that both $J$ and $K$ are principal. Thus, any ideal that divides a principal ideal of a domain $R$ with trivial Picard group is itself principal.
\epf

\bc \label{coro2.5}
A molecular domain with trivial Picard group is atomic.
\ec

\bpf
Let $R$ be a molecular domain with trivial Picard group and let $a$ be a nonzero nonunit of $R$. Since $R$ is molecular, there exists a molecularization $(a) = I_1I_2 \, \cdots \, I_n$. By Proposition \ref{prop2} above, each $I_j$ is principal and, moreover, is generated by an atom $a_j$ of $R$. Hence, $a = ua_1a_2 \, \cdots \, a_n$ for some unit $u$, and so $a$ can be written as a product of atoms, as desired.
\epf

Recall that a \textit{multiplication ideal} of a ring $R$ is an ideal $I$ of $R$ satisfying the property that for any ideal $J \subseteq I$ of $R$, there exists an ideal $K$ of $R$ for which $J = IK$. Note that an invertible ideal is a multiplication ideal. Proposition \ref{prop0} below, while elementary, highlights the intuitive idea that ``multiplication ideal" and ``molecule" are dynamically dissimilar in terms of the degree of ideal factorization involved.

\bp \label{prop0}
Let $R$ be a domain and $I$ a molecule of $R$. If $J$ is a proper ideal of $R$ for which $I \subsetneq J$, then $J$ cannot be a multiplication ideal of $R$.
\ep

\bpf Deny. Then there exists a necessarily proper ideal $K$ of $R$ for which $I = JK$, contradicting the fact that $I$ is a molecule of $R$. \epf

Certainly, the class of Dedekind domains reveals that sometimes the only ideals of a domain that are molecules are the maximal ideals themselves. However, maximal ideals need not be molecules, for instance, as in a valuation domain with a non-principal maximal ideal (in fact, such a domain has no molecules at all).

It is natural then to seek out domains where a cancellation-type property for ideals holds, as such contexts can give rise to the existence of molecules from which an exploration of associated factorizations can begin. We find an abundance of fruitfulness in considering domains that have what we deem ``unit-cancellation for ideals". Motivated by the ``pr\'{e}simplifiable condition" for commutative rings with identity ($xy = x \Rightarrow x = 0 \textrm{ or } y \textrm{ is a unit}$), we say that an ideal $I$ of the domain $R$ is {\em unit-cancellative} if for each ideal $J$ of $R$ with $I = IJ$, it must be the case that $J = R$. By extension, we say the domain $R$ has {\em unit-cancellation for ideals} if every nonzero ideal of $R$ is unit-cancellative. This notion is equivalent to the monoid of nonzero ideals of $R$ being ``unit-cancellative" in the terminology of \cite{FGKT}.

Propositions \ref{prop4.51} and \ref{prop4.5} below provide for a wealth of domains that have unit-cancellation for ideals. Recall that if $S$ is a ring, $R$ a subring of $S$, and $I$ an ideal of $R$, then we say that $I$ {\em survives} in $S$ if $IS$ is a proper ideal of $S$ and, moreover, we say that $R \subseteq S$ is {\em survival extension} if every proper ideal of $R$ survives in $S$. Clearly, every integral extension is a survival extension, owing to the lying-over theorem. Also, if $S$ is quasisemilocal and $U(S) \cap R = U(R)$, then an application of prime avoidance reveals that $R \subseteq S$ is a survival extension.

\bp \label{prop4.51}
Let $\Omega$ be a nonempty set of domains that have unit-cancellation for ideals and let $R$ be a domain that is a subring of each $S \in \Omega$. If for each proper ideal $I$ of $R$ there is some $S \in \Omega$ such that $I$ survives in $S$, then $R$ has unit-cancellation for ideals.
\ep

\bpf
Let $I$ be a nonzero ideal of $R$ and $J$ an ideal of $R$ for which $I = IJ$. Assume that $J$ is proper. Then $JS \neq S$ for some $S \in \Omega$. On the other hand, $IS = ISJS$, and hence $JS = S$, a contradiction.
\epf

\bc \label{coro4.52}
If $R \subseteq S$ is a survival extension of domains and $S$ has unit-cancellation for ideals, then $R$ has unit-cancellation for ideals.
\ec

\bp \label{prop4.5}
Let $R$ be a domain. If $R$ satisfies the conclusion of the Krull intersection theorem (that is, $\cap_{n=0}^{\infty} M^n = 0$ for all maximal ideals $M$ of $R$), then $R$ has unit-cancellation for ideals. In particular, if $R$ is Noetherian or the integral closure $\overline{R}$ in some field extension of the quotient field of $R$ is Noetherian, then $R$ has unit-cancellation for ideals.
\ep

\bpf
Let $R$ be a domain that satisfies the conclusion of the Krull intersection theorem. Let $I$ be a nonzero ideal of $R$ and $J$ an ideal of $R$ for which $I = IJ$. Observe that $I = IJ^n$ for every whole number $n$, and so $I \subseteq \cap_{n=0}^{\infty} J^n$. We conclude that $J = R$. Therefore, $R$ has unit-cancellation for ideals.

Since it is well-known that every Noetherian domain satisfies the conclusion of the Krull intersection theorem, it now follows that if $R$ is a Noetherian domain, then $R$ has unit-cancellation for ideals (an alternative justification of the fact that Noetherian domains have unit-cancellation for ideals is found by observing that any nonzero finitely generated ideal of a domain is unit-cancellative, an application of Nakayama's lemma being all that is necessary to show this). Suppose then that the integral closure $\overline{R}$ in some field extension of the quotient field of $R$ is Noetherian and let $M$ be a maximal ideal of $R$. By the lying-over theorem, there is some maximal ideal $\mathcal{M}$ of $\overline{R}$ such that $\mathcal{M} \cap R = M$. Since $\overline{R}$ satisfies the conclusion of the Krull intersection theorem, we have that $\cap_{n=0}^{\infty} M^n \subseteq \cap_{n=0}^{\infty} \mathcal{M}^n = 0$, and so $\cap_{n=0}^{\infty} M^n = 0$. Therefore, $R$ itself must satisfy the conclusion of the Krull intersection theorem from which we deduce that $R$ itself has unit-cancellation for ideals. The proof is thus complete.
\epf

It should be noted that a polynomial ring in infinitely many indeterminates over a field has unit-cancellation for ideals, but the ring is not Noetherian nor is any integral extension of the ring Noetherian.

We now give a proposition (Proposition \ref{prop6.1}) that characterizes molecules in terms of the unit-cancellation property. As a consequence, domains with unit-cancellation for ideals adequately address an issue mentioned in the discussion just prior to Proposition \ref{prop4.5}.

\bp \label{prop6.1}
Let $R$ be a domain, $I$ a proper ideal of $R$, $P$ a prime ideal of $R$, and $M$ a maximal ideal of $R$.

(a) $I$ is a molecule of $R$ if and only if (1) $I$ is unit-cancellative and (2) for all ideals $J$ and $K$ of $R$ for which $I = JK$, it follows that either $J = I$ or $K=I$.

(b) $P$ is a molecule if and only if $P$ is unit-cancellative.

(c) (cf. \cite[Corollary 2.3]{HL1}) $M$ is a molecule if and only if $M$ is not idempotent.
\ep

\bpf
(a) Trivially, if $I$ is a molecule of $R$, then $I$ is unit-cancellative and for all ideals $J$ and $K$ of $R$ such that $I = JK$, it must be the case that either $J = I$ or $K = I$. Conversely, let $J$ and $K$ be ideals of $R$ such that $I = JK$. Then either $J = I$ or $K = I$ by hypothesis. Without loss of generality, assume that $J = I$. Then $I = IK$, and so $K = R$ since $I$ is unit-cancellative. Therefore, $I$ is a molecule of $R$.

(b) Clearly, if $J$ and $K$ are ideals of $R$ such that $P = JK$, then either $J = P$ or $K = P$. Therefore, the statement is an immediate consequence of (a).

(c) If $M$ is unit-cancellative, then $M$ is clearly not idempotent. If $M$ is not idempotent and $J$ is an ideal of $R$ such that $M = MJ$, then $M$ is properly contained in $J$, whence $J = R$, and so $M$ is unit-cancellative. The statement now follows from (b).
\epf

Proposition \ref{prop7} below reveals that molecules share a property of strongly irreducible ideals in regards to products of pairwise comaximal ideals. It is this consideration of products of a certain type of ideal that also gives rise to a large class of domains where every molecule is a primary ideal (cf. \cite[Example 2.14]{HL1}), formalized in Proposition \ref{prop7.1}.

\bp \label{prop7}
Let $R$ be a domain and $I$ a molecule of $R$. Let $J_1,J_2, \, \ldots \, ,J_n$ be pairwise comaximal ideals of $R$ such that $J_1J_2 \, \cdots \, J_n \subseteq I$. Then $J_i \subseteq I$ for some $i =1,2, \, \ldots \, ,n$.
\ep

\bpf
Note that if $J_1,J_2, \, \ldots \, ,J_n$ are pairwise comaximal ideals of $R$, then $\sum \widehat{J_1J_2 \, \cdots \, J_n} = R$. As such $I = (I + J_1)(I+J_2) \, \cdots \, (I+J_n)$. However, since $I$ is a molecule, $I+J_i = I$ for some $i =1,2, \, \ldots \, ,n$. This means that $J_i \subseteq I$, as desired.
\epf

\bp \label{prop7.1}
Let $R$ be a Laskerian domain such that every non-maximal prime ideal of $R$ is a multiplication ideal of $R$. Then every molecule of $R$ is a primary ideal of $R$.
\ep

\bpf
By \cite[Theorem 10]{AM0}, every proper ideal of $R$ is a finite product of primary ideals of $R$. It is obvious then that every molecule of $R$ is a primary ideal of $R$.
\epf

We can also specialize to sufficient conditions on a molecule itself that guarantee that the molecule is primary, as given in Proposition \ref{prop7.11} and Corollary \ref{coro7.12} below. This has the upshot of providing an alternative means of obtaining a result on molecules of a one-dimensional Noetherian domain (Corollary \ref{coro8}) that is also a consequence of Proposition \ref{prop7.1}.

\bp \label{prop7.11}
Let $R$ be a domain and $I$ a molecule of $R$ such that $R/I$ is a quasisemilocal zero-dimensional ring whose nilradical is nilpotent. Then $I$ is a primary ideal of $R$.
\ep

\bpf Since $R/I$ is quasisemilocal and zero-dimensional, we have that the set of all prime ideals of $R$ that contain $I$ is finite and consists only of maximal ideals of $R$. Let $\{M_1,M_2, \, \ldots \, ,M_n\}$ be the set of all prime ideals of $R$ that contain $I$. Since the nilradical of $R/I$, which is $\sqrt{I}/I$, is nilpotent, it follows that $(\sqrt{I})^m \subseteq I$ for some $m \in \mathbb{N}$. However, this means that $\prod_{i=1}^n M_i^m \subseteq I$, and so $M_j^m \subseteq I$ for some $j=1,2, \, \ldots \, ,n$ by Proposition \ref{prop7}. We conclude then that $I$ is $M_j$-primary.
\epf

\bc \label{coro7.12}
Let $R$ be a domain and $I$ a molecule of $R$ such that $R/I$ is Artinian. Then $I$ is a primary ideal of $R$.
\ec

\bc \label{coro8}
Let $R$ be a one-dimensional Noetherian domain. Then every molecule of $R$ is a primary ideal of $R$.
\ec

Proposition \ref{prop9} below provides one of the key reasons for opting for the ``molecule" terminology over the previous ``factorable" language. In particular, noting Butts' characterization of a domain with unique factorization of ideals into nonfactorable ideals as a Dedekind domain \cite[Theorem]{HSB}, one should recognize the parallel of the statement of Proposition \ref{prop9} with the well-known characterization of UFD's as atomic domains where every atom is prime.

\bp \label{prop9}
Let $R$ be a domain. Then $R$ is a molecular domain for which every molecule is prime if and only if $R$ is a Dedekind domain.
\ep

\bpf
The fact that if $R$ is a Dedekind domain, then $R$ is a molecular domain for which every molecule is prime was established by Butts \cite{HSB}. Conversely, if $R$ is a molecular domain for which every molecule is prime, then every ideal of $R$ is a product of prime ideals of $R$, and it is well-known that such a domain must be Dedekind.
\epf

\section{FINITE MOLECULARIZATION DOMAINS} \label{S2.1}

We now come to the main concept of this paper, the notion of a ``finite molecularization domain", an idea designed to be the ideal-theoretic analogue of a ``finite factorization domain".

\bd \label{defn0}
A {\em finite molecularization domain}, or FMD, is a molecular domain $R$ with the property that every nonzero proper ideal of $R$ has only a finite number of molecularizations--that is, every nonzero proper ideal of $R$ has only a finite number of factorizations as a product of molecules. Moreover, if $I$ is an ideal of $R$ such that $I = J_1J_2 \, \cdots \, J_m$ and $I = K_1K_2 \, \cdots \, K_n$, where each $J_i$ and $K_j$ is a molecule of $R$, then these factorizations will be regarded as the same if $m=n$ and there is a permutation $\sigma \in S_n$ for which $J_i = K_{\sigma(i)}$ for $i = 1,2, \, \ldots \, ,n$.
\ed

In light of the obvious connection that the notion of ``finite molecularization domain" has with ``finite factorization domain", it is natural to ask if the property of being an FMD can be characterized, at least in part, by the finiteness of the number of (non-associated) molecular divisors (see \cite[p. 2]{AAZ1}). Corollary \ref{coro1.01} below reveals an affirmative answer to this question. In fact, the characterization presented through Corollary \ref{coro1.01} is powerful enough to spawn several corollaries of its own, two of which (Corollaries \ref{coro1.2} and \ref{coro1.3}) establish that every FSD is an FMD and every FMD is an FFD.

We first lay the groundwork for Corollary \ref{coro1.01} through Theorem \ref{theo1}, which addresses the situation at the level of the ideals themselves.

\bt \label{theo1}
Let $R$ be a molecular domain and $I$ a unit-cancellative ideal of $R$ that is divisible by only finitely many molecules of $R$. Then $I$ has only a finite number of molecularizations (that is, $I$ has only a finite number of factorizations as a product of molecules of $R$) and $I$ is divisible by only finitely many ideals of $R$.
\et

\bpf
Let $R$ be a molecular domain and $I$ a unit-cancellative ideal of $R$ that is divisible by only finitely many molecules of $R$. Without loss of generality, we may assume that $I$ is a nonzero proper ideal of $R$. Let $\{J_1, J_2, \, \ldots \, ,J_n\}$ be the set of all molecules that divide $I$. Put $E = \{(\alpha_1, \alpha_2, \, \ldots \, ,\alpha_n) \in \mathbb{W}^n \mid I = \prod_{i=1}^n J_i^{\alpha_i}\}$. Then $\mbox{{\rm Min}}(E)$ is finite by Dickson's lemma.

Next, we show that $\mbox{{\rm Min}}(E) = E$. Let $\alpha = (\alpha_1, \alpha_2, \, \ldots \, ,\alpha_n) \in E$. By Dickson's lemma, there is some $\beta = (\beta_1, \beta_2, \, \ldots \, ,\beta_n) \in \mbox{{\rm Min}}(E)$ such that $\beta \leq \alpha$ (that is, $\beta_i \leq \alpha_i$ for each $i = 1,2, \, \ldots \, ,n$). Thus, $I = \prod_{i=1}^n J_i^{\alpha_i} = \prod_{i=1}^n J_i^{\beta_i} \prod_{i=1}^n J_i^{\alpha_i - \beta_i} = I \prod_{i=1}^n J_i^{\alpha_i - \beta_i}$ and hence $\prod_{i=1}^n J_i^{\alpha_i - \beta_i} = R$. As such, $\alpha = \beta$, and so $\mbox{{\rm Min}}(E) = E$. It follows that $E$ is finite, and thus $I$ has only a finite number of factorizations as a product of molecules of $R$. Furthermore, observe that the set of ideals of $R$ that divide $I$ is given by $\{\prod_{i=1}^n J_i^{\gamma_i} \mid \gamma = (\gamma_1, \gamma_2, \, \ldots \, ,\gamma_n) \in \mathbb{W}^n \textrm{ and } \gamma \leq \alpha \textrm{ for some } \alpha \in E\}$, which is clearly finite. The proof is thus complete.
\epf

\bc \label{coro1.01}
Let $R$ be a domain. The following are equivalent:

(1) $R$ is an FMD;

(2) $R$ is a molecular domain, $R$ has unit-cancellation for ideals, and every nonzero ideal of $R$ is divisible by only finitely many molecules of $R$;

(3) $R$ has unit-cancellation for ideals and every nonzero ideal of $R$ is divisible by only finitely many ideals of $R$.
\ec

\bpf
(1) $\Rightarrow$ (2): Obviously, $R$ is molecular. Let $I$ be a nonzero proper ideal of $R$. Clearly, $I$ is divisible by only finitely many molecules of $R$. Let $\{J_1, J_2, \, \ldots \, ,J_n\}$ be the set of all molecules that divide $I$. Since $R$ is an FMD, $E = \{(\gamma_1, \gamma_2, \, \ldots \, ,\gamma_n) \in \mathbb{W}^n \mid I = \prod_{i=1}^n J_i^{\gamma_i}\}$ is finite. In particular, there is some $\alpha = (\alpha_1, \alpha_2, \, \ldots \, ,\alpha_n) \in \mathbb{W}^n$ such that $\gamma \leq \alpha$ for all $\gamma \in E$. Let $J$ be an ideal of $R$ such that $I = IJ$. Then there are $\beta = (\beta_1, \beta_2, \, \ldots \, ,\beta_n) \in \mathbb{W}^n$ and $\gamma =(\gamma_1, \gamma_2, \, \ldots \, ,\gamma_n) \in \mathbb{W}^n$ for which $I = \prod_{i=1}^n J_i^{\beta_i}$ and $J = \prod_{i=1}^n J_i^{\gamma_i}$. However, if $m \in \mathbb{W}$, then $I = IJ^m$, whence $I = \prod_{i=1}^n J_i^{\beta_i + m\gamma_i}$. But this means that $\beta + m\gamma \leq \alpha$ for every $m \in \mathbb{W}$, and so $\gamma_i = 0$ for each $i = 1,2, \, \ldots \, ,n$. Therefore, $J = R$, and we conclude that $I$ is unit-cancellative.

(2) $\Rightarrow$ (1) and (3): This follows immediately from Theorem \ref{theo1} above.

(3) $\Rightarrow$ (2): It is sufficient to show that $R$ is a molecular domain. Let $I$ be a nonzero proper ideal of $R$. Assume to the contrary that $I$ is not a finite product of molecules of $R$. Clearly, there is a proper ideal $J$ of $R$ that divides $I$ such that $J$ is maximal amongst the proper ideals of $R$ that divide $I$ and that are not a finite product of molecules of $R$. Since $J$ itself cannot be a molecule of $R$, it follows by Proposition \ref{prop6.1} that there are proper ideals $A$ and $B$ of $R$ such that $J \subsetneq A,B$ and $J = AB$. But then $A$ and $B$ are each finite products of molecules of $R$, whence $J$ must be too, a contradiction. This completes the proof.
\epf

\bc \label{coro1.3}
Every FMD is an FFD.
\ec

\bpf
From Corollary \ref{coro1.01}, it follows that every nonzero principal ideal of the FMD $R$ is contained in only finitely many principal ideals of $R$, and thus $R$ is an FFD by \cite[Theorem 1]{AMul}.
\epf

\bc \label{coro1.1}
Let $R$ be a Noetherian domain. Then $R$ is an FMD if and only if every nonzero ideal of $R$ is divisible by only finitely many molecules of $R$.
\ec

\bc \label{coro1.2}
Every FSD is an FMD.
\ec

\bpf
Let $R$ be an FSD. Then $R$ is Noetherian \cite[Proposition 2.1]{HL0} and every nonzero ideal of $R$ is contained in (and thus divisible by) only finitely many molecules of $R$. Therefore, $R$ is an FMD by Corollary \ref{coro1.1}.
\epf

In light of Corollaries \ref{coro1.1} and \ref{coro1.2} above, it should be noted that while every Noetherian domain is molecular (see \cite[Theorem 2.13]{HL1}), not every Noetherian domain is an FMD, as revealed in Theorems \ref{theo10} and \ref{theo16}. In fact, Theorem \ref{theo16} provides for the existence of one-dimensional local domains that are not FMD's. Nonetheless, since it is well-known that FFD's satisfy ACCP, Corollary \ref{coro1.3} gives that FMD's also satisfy ACCP.

It also should be pointed out that the converse of the implication in Corollary \ref{coro1.3} is false, in general. Since \cite[Theorem 5]{AKP} implies that a Pr\"{u}fer FMD is a Dedekind domain, any one-dimensional Pr\"{u}fer FFD that is not Dedekind (see, amongst others, \cite[Example 2]{G}) cannot be an FMD.

Now, given that the property of unit-cancellation for ideals and the property that every nonzero ideal is divisible by only finitely many ideals were used to characterize FMD's in Corollary \ref{coro1.01}, we offer the following result in the spirit of the domain extension considerations of Proposition \ref{prop4.51} and Corollary \ref{coro4.52}.

\bp \label{prop19.98}
Let $S$ be a domain and $R$ a subring of $S$ such that the conductor ideal $(R : S) \neq 0$. If every nonzero ideal of $R$ is divisible by only finitely many ideals of $R$, then every nonzero ideal of $S$ is divisible by only finitely many ideals of $S$.
\ep

\bpf
Let $I$ be a nonzero ideal of $S$ and let $x \in (R : S)$ be nonzero. Observe that if $J$ and $L$ are ideals of $S$ such that $I = JL$, then $x^2I = xJxL$ and $x^2I$, $xJ$, and $xL$ are ideals of $R$. Let $\mathcal{I}$ be the set of ideals of $S$ that divide $I$ and let $\mathcal{J}$ be the set of ideals of $R$ that divide $x^2I$. Let $f$:$\,\,\mathcal{I} \to \mathcal{J}$ be defined by $f(J) = xJ$. Since $f$ is a well-defined injective map, we must have that $\mathcal{I}$ is finite, as desired.
\epf

We now transition to local-global considerations with regards to the property of being an FMD. Our main result along these lines, Theorem \ref{theo20}, shows that the property of being an FMD is stable at localizations of height-one maximal ideals (as contrasted with the property of being an FFD which is not, in general, stable under the formation of localizations at height-one maximal ideals; see \cite[Example 2]{G} and \cite[Example 5.4]{AAZ1}). We first provide a lemma regarding some special properties of ideals of localizations at height-one maximal ideals.

\bl \label{lemm19.99}
Let $R$ be a domain, $M$ a height-one maximal ideal of $R$, and $I$ an ideal of $R_M$.

(a) If $C$ and $D$ are ideals of $R_M$ such that $I = CD$, then $I \cap R = (C \cap R)(D \cap R)$.

(b) If $I$ is principal, then $I \cap R$ is locally principal.

(c) $I$ is the only ideal of $R_M$ whose contraction to $R$ is $I \cap R$.
\el

\bpf (a) Let $C$ and $D$ be ideals of $R_M$ such that $I = CD$. Without loss of generality, we may assume that $I$ is a nonzero proper ideal of $R_M$. Put $A = C \cap R$ and $B = D \cap R$. Since $I$ is a primary ideal of $R_M$, we have that $I \cap R$ is an $M$-primary ideal of $R$. Moreover, $I \cap R \subseteq A \cap B$, and thus $M = \sqrt{I \cap R} \subseteq \sqrt{A} \cap \sqrt{B} = \sqrt{AB}$. If $\sqrt{AB} = R$, then $A = B = R$, whence $I = R_M$, a contradiction. Therefore, $\sqrt{AB} = M$, and so $AB$ is $M$-primary. Since $C = A_M$ and $D = B_M$, we conclude that $I \cap R = (AB)_M \cap R = AB = (C \cap R)(D \cap R)$.

(b) Let $I$ be principal and nonzero, and let $N$ be a maximal ideal of $R$. Note that $M \subseteq \sqrt{I \cap R}$. If $N \neq M$, then $I \cap R$ is not contained in $N$, and hence $(I \cap R)_N = R_N$ is principal. On the other hand, if $N = M$, then $(I \cap R)_N = I$ is principal.

(c) Let $J$ be an ideal of $R_M$ such that $J \cap R = I \cap R$. Then $J = (J \cap R)_M = (I \cap R)_M = I$.
\epf

\bt \label{theo20}
Let $R$ be a domain and $M$ a height-one maximal ideal of $R$.

(a) If $R$ has unit-cancellation for ideals, then $R_M$ has unit-cancellation for ideals.

(b) If every nonzero ideal of $R$ is divisible by only finitely many ideals of $R$, then every nonzero ideal of $R_M$ is divisible by only finitely many ideals of $R_M$.

(c) If $R$ is an FMD, then $R_M$ is an FMD.
\et

\bpf
(a) Let $I$ be a nonzero ideal of $R_M$ and $J$ an ideal of $R_M$ such that $I = IJ$. By Lemma \ref{lemm19.99}(a), we have that $I \cap R = (I \cap R)(J \cap R)$. Clearly, $I \cap R$ is nonzero, and thus $I \cap R$ is a unit-cancellative ideal of $R$. It follows that $J \cap R = R$, and so $J = R_M$.

(b) Let $I$ be a nonzero ideal of $R_M$. Let $\mathcal{I}$ be the set of ideals of $R_M$ that divide $I$ and $\mathcal{J}$ the set of ideals of $R$ that divide $I \cap R$. Let $f$:$\,\,\mathcal{I} \to \mathcal{J}$ be given by $f(L) = L \cap R$. Then $f$ is a well-defined injective map by Lemma \ref{lemm19.99}. Since $\mathcal{J}$ is finite, we must have that $\mathcal{I}$ is finite.

(c) This is an immediate consequence of (a) and (b) above and Corollary \ref{coro1.01}.
\epf

Note that a certain converse of Theorem \ref{theo20}(a) is true, in the sense that if $R_M$ has unit-cancellation for ideals for every maximal ideal $M$ of $R$, then $R$ has unit-cancellation for ideals, thanks to Proposition \ref{prop4.51}.

While we can certainly deduce from Theorem \ref{theo20} that if $R$ is an FMD and $M$ is a height-one maximal ideal of $R$, then $R_M$ is an FFD (see Corollary \ref{coro1.3}), Proposition \ref{prop20.02} below shows that this conclusion can actually be obtained under a slightly more general hypothesis.

\bp \label{prop20.02}
Let $R$ be a domain for which each nonzero locally principal ideal is divisible by only finitely many locally principal ideals. If $M$ is a height-one maximal ideal of $R$, then $R_M$ is an FFD.
\ep

\bpf
Let $I$ be a nonzero principal ideal of $R_M$. Let $\mathcal{I}$ be the set of principal ideals of $R_M$ that divide $I$ and $\mathcal{J}$ the set of locally principal ideals of $R$ that divide $I \cap R$. Let $f$:$\,\,\mathcal{I} \to \mathcal{J}$ be given by $f(L) = L \cap R$. Then $f$ is a well-defined injective map by Lemma \ref{lemm19.99}. Since $\mathcal{J}$ is finite, we must have that $\mathcal{I}$ is finite. Since principal ideals are multiplication ideals, there are only finitely many principal ideals of $R_M$ that contain $I$. Therefore, $R_M$ is an FFD by \cite[Theorem 1]{AMul}.
\epf

We conclude this section making contact again with the FSD property. In particular, after presenting a novel characterization through Proposition \ref{prop22} below of when the property of being an FSD globalizes, we are able to present in Theorem \ref{theo23} a characterization of FSD's in terms of the FMD property. As a consequence of this theorem (Corollary \ref{coro23.01}), we are able to answer a conjecture presented in \cite{HL0} regarding sufficient conditions guaranteeing the equivalence of the FSD and FFD properties.

\bp \label{prop22}
Let $R$ be a domain. Then $R$ is an FSD if and only if $R$ is of finite character and locally an FSD.
\ep

\bpf
Let $R$ be an FSD. Clearly, $R$ is of finite character. It also follows by \cite[Theorem 3.5]{HL0} that $R_M$ is an FSD for each maximal ideal $M$ of $R$.

Conversely, let $R$ be of finite character and locally an FSD. Let $I$ be a nonzero proper ideal of $R$. Let $M_1, M_2, \, \ldots \, ,M_n$ be the maximal ideals of $R$ that contain $I$. For each $i =1,2, \, \ldots \, ,n$, let $\mathcal{J}_i$ be the set of ideals of $R_{M_i}$ that contain $I_{M_i}$ and let $\mathcal{J}$ be the set of ideals of $R$ that contain $I$. If $J \in \mathcal{J}$ and $\mbox{{\rm Max}}(R)$ is the set of all maximal ideals of $R$, then $J = R \cap \left(\cap_{M \in \mbox{{\rm Max}}(R)} J_M \right) = R \cap \left(\cap_{i=1}^n J_{M_i}\right)$. Therefore, $f$:$\,\,\prod_{i=1}^n \mathcal{J}_i \to \mathcal{J}$ defined by $f((J_i)_{i=1}^n) = R \cap \left(\cap_{i=1}^n J_i\right)$ is a well-defined surjective map. However, since $\mathcal{J}_i$ is finite for each $i =1,2, \, \ldots \, ,n$, it must be the case that $\mathcal{J}$ is finite, as desired.
\epf

We pause briefly to expand upon the value of the ``finite character" hypothesis utilized in Proposition \ref{prop22} above by recognizing its worth for FMD considerations through Proposition \ref{prop22.1} below.

\bp \label{prop22.1}
Let $R$ be a domain of finite character.

(a) If for each maximal ideal $M$ of $R$ we have that every nonzero ideal of $R_M$ is divisible by only finitely many ideals of $R_M$, then every nonzero ideal of $R$ is divisible by only finitely many ideals of $R$.

(b) If $R$ is locally an FMD, then $R$ is an FMD.
\ep

\bpf
(a) Let $I$ be a nonzero ideal of $R$. Let $\mathcal{I}$ be the set of all ideals of $R$ that divide $I$, and for each maximal ideal $M$ of $R$, let $\mathcal{I}_M$ be the set of all ideals of $R_M$ that divide $I_M$. Let $\mbox{{\rm Max}}(R)$ be the set of all maximal ideals of $R$, and let $f$:$\,\,\mathcal{I} \to \prod_{M \in \mbox{{\rm Max}}(R), \, I \subseteq M} \mathcal{I}_M$ be given by $f(L) = \left(L_M\right)_{M \in \mbox{{\rm Max}}(R), \, I \subseteq M}$. Clearly, $f$ is a well-defined map. Let $J,L \in \mathcal{I}$ be such that $f(J) = f(L)$ and let $M \in \mbox{{\rm Max}}(R)$. If $I \subseteq M$, then we clearly have that $J_M = L_M$. If $I \nsubseteq M$, then $J \nsubseteq M$ and $L \nsubseteq M$, and so $J_M = R_M = L_M$. We deduce that $J = L$, and thus $f$ is injective. It follows that $\mathcal{I}$ is finite, as desired.

(b) This follows from (a) above, Proposition \ref{prop4.51}, and Corollary \ref{coro1.01}.
\epf

We note that the ``finite character" hypothesis cannot be dropped from Proposition \ref{prop22.1} above by considering an almost Dedekind domain $R$ that is not a Dedekind domain. For such a domain $R$ is locally an FSD, hence locally an FMD, but, since $R$ is a Pr\"{u}fer domain that is not a Dedekind domain, it cannot be an FMD (see \cite[Theorem 5]{AKP}).

We now give the promised characterization of FSD's in terms of the FMD property.

\bt \label{theo23}
Let $R$ be a domain. The following are equivalent:

(1) $R$ is an FSD;

(2) $R$ is a Noetherian FMD and $\textrm{dim}(R) \leq 1$;

(3) $R$ is Noetherian, $\textrm{dim}(R) \leq 1$, and every invertible ideal of $R$ is contained in only finitely many invertible ideals of $R$.
\et

\bpf
(1) $\Rightarrow$ (2): Let $R$ be an FSD. By \cite[Proposition 2.1]{HL0}, we have that $R$ is Noetherian and $\textrm{dim}(R) \leq 1$. It follows from Corollary \ref{coro1.2} that $R$ is an FMD, as well.

(2) $\Rightarrow$ (3): Let $I$ be an invertible ideal of $R$. By Corollary \ref{coro1.01}, $I$ is divisible by only finitely many (invertible) ideals of $R$. The assertion now follows since invertible ideals are multiplication ideals.

(3) $\Rightarrow$ (1): Without loss of generality, we may assume that $R$ is not a field. Since $R$ is Noetherian, we have that every nonzero locally principal ideal of $R$ is invertible. If $M$ is a maximal ideal of $R$, then $R_M$ is Noetherian, and hence $R_M$ is an FSD by combining Proposition \ref{prop20.02} with \cite[Theorem 2.5]{HL0}. Clearly, $R$ is of finite character, and so $R$ must be an FSD by Proposition \ref{prop22}.
\epf

\bc \label{coro23.01}
Let $R$ be a Noetherian one-dimensional domain whose Picard group is trivial. Then $R$ is an FSD if and only if $R$ is an FFD.
\ec

As a result of Corollary \ref{coro23.01} above, the first conjecture expressed in the remark immediately following \cite[Theorem 2.5]{HL0} is proved. In particular, if $R$ is a one-dimensional (Noetherian) semilocal domain, then $R$ is an FSD if and only if $R$ is an FFD. Moreover, this equivalence breaks down if the ``semilocal" hypothesis is removed, as evidenced by the fact that $k[X^2,X^3]$, where $k$ is an infinite field, is a one-dimensional Noetherian FFD that is not an FSD (see Theorem \ref{theo10} and \cite[Theorem 3]{AMul}).

Let $R$ be a domain. Recall that $R$ is called {\em seminormal} if whenever $x$ is an element of the quotient field of $R$ such that $x^2,x^3 \in R$, it must be the case that $x \in R$, and the {\em seminormalization} of $R$ is the intersection of all seminormal overrings of $R$. We will denote the complete integral closure of $R$ by $\widehat{R}$. As in \cite{Gil}, we say that $R$ is a {\em G}-domain if the intersection of all nonzero prime ideals of $R$ is nonzero.

Our final result of this section, Theorem \ref{theo24.01}, improves upon the characterization of FSD's offered in Theorem \ref{theo23} under the additional hypothesis that $(R : S) \neq 0$, where $S$ is the seminormalization of $R$.

\bt \label{theo24.01}
Let $R$ be a domain and $S$ the seminormalization of $R$. If the conductor ideal $(R : S) \neq 0$, then the following are equivalent:

(1) $R$ is an FSD;

(2) $R$ is an FMD and $\overline{R}$ is a Dedekind domain;

(3) $R$ is a Mori domain, $\textrm{dim}(R) \leq 1$, and every invertible ideal of $R$ is contained in only finitely many invertible ideals of $R$;

(4) $\textrm{dim}(R) \leq 1$, $\widehat{R}$ is a Krull domain, and every nonzero locally principal ideal of $R$ is divisible by only finitely many locally principal ideals of $R$.
\et

\bpf
We begin by claiming that if $(R : S) \neq 0$, where $S$ is the seminormalization of $R$, then $(R_P : \widehat{R_P}) \neq 0$ for each height-one prime ideal $P$ of the domain $R$. For let $P$ be a height-one prime ideal of the domain $R$, let $S$ be the seminormalization of $R$, and suppose that $(R : S) \neq 0$. Note that $(R_P : S_P) \supseteq (R : S)$ and $S_P \subseteq \overline{R}_P = \overline{R_P} \subseteq \widehat{R_P} \subseteq \widehat{S_P}$. Since $R_P$ is quasilocal and one-dimensional, it must be a {\em G}-domain. We have that $S_P$ is an overring of $R_P$, and thus $S_P$ is a {\em G}-domain, as well. Moreover, $S_P$ is seminormal, and so it follows by \cite[Proposition 4.8.2]{GHHK} that $(S_P : \widehat{R_P}) \supseteq (S_P : \widehat{S_P}) \neq 0$. We conclude that $(R_P : \widehat{R_P}) \supseteq (R_P : S_P)(S_P : \widehat{R_P}) \neq 0$, and the claim is thus proved.

(1) $\Rightarrow$ (2) and (3): This follows immediately from Theorem \ref{theo23} and the Krull-Akizuki theorem.

(2) $\Rightarrow$ (4): It follows from Corollary \ref{coro1.01} that every nonzero ideal of $R$ is divisible by only finitely many ideals of $R$, and so the corresponding statement about locally principal ideals of $R$ is clearly true. Observe that $\textrm{dim}(R) = \textrm{dim}(\overline{R}) \leq 1$. Since $\overline{R}$ is Noetherian, we have that $\overline{R} \subseteq \widehat{R} \subseteq \widehat{\overline{R}} = \overline{R}$. Consequently, $\widehat{R}$ is a Dedekind domain, and hence a Krull domain.

(3) $\Rightarrow$ (1): Without loss of generality, we may assume that $R$ is not a field. Since $R$ is a one-dimensional Mori domain, it is clear that $R$ is of finite character.

Next, we show that every nonzero locally principal ideal of $R$ is invertible. Let $I$ be a nonzero locally principal of $R$ and $N$ a maximal ideal of $R$. Since $R$ is a Mori domain, there is some nonzero finitely generated ideal $J$ of $R$ such that $J \subseteq I$ and $I^{-1} = J^{-1}$. Therefore, $(I^{-1})_N = (J^{-1})_N = J_N^{-1} \supseteq I_N^{-1} \supseteq (I^{-1})_N$. We deduce that $(II^{-1})_N = I_NI_N^{-1} = R_N$, and hence $II^{-1} = R$.

Now, let $M$ be a maximal ideal of $R$. By Proposition \ref{prop22}, it suffices to show that $R_M$ is an FSD. Observe that $R_M$ is a Mori domain. The claim proved at the beginning of this proof implies that $(R_M : \widehat{R_M}) \neq 0$. By Proposition \ref{prop20.02}, $R_M$ is an FFD, and thus $U(\widehat{R_M})/U(R_M)$ is finite by \cite[Theorem 4]{AMul}. It follows by \cite[Theorem 4.2]{R} that $R_M$ is Noetherian, and thus $R_M$ is an FSD by \cite[Proposition 2.16]{HL0}.

(4) $\Rightarrow$ (1): Without loss of generality, we may assume that $R$ is not a field. Let $M$ be a maximal ideal of $R$. By Proposition \ref{prop20.02}, $R_M$ is an FFD and, by the claim proved at the beginning of this proof, we have that $(R_M : \widehat{R_M}) \neq 0$. Observe that $\widehat{R}_M$ is a Krull domain and $\widehat{R}_M \subseteq \widehat{R_M}$. Since $R_M \subseteq \widehat{R}_M$ and $\widehat{R}_M$ is completely integrally closed, we have that $\widehat{R_M} \subseteq \widehat{R}_M$. Therefore, $\widehat{R_M} = \widehat{R}_M$ is a Krull domain. Since $R_M$ is a {\em G}-domain, $\widehat{R_M}$ is a {\em G}-domain, and thus $\widehat{R_M}$ is a semilocal PID (cf. \cite[Lemma 3.1]{R}). We conclude by \cite[Lemma 5.2(1) and Proposition 5.6(2)]{GR} that $R_M$ is Noetherian, and thus $R_M$ is an FSD by \cite[Proposition 2.16]{HL0}.

Since $(R_M : \widehat{R_M}) \neq 0$, we have that $\widehat{R_M}$ is not a field. Let $Q$ be a maximal ideal of $\widehat{R_M}$. Then $Q \cap R_M = M_M$. Since $\widehat{R} \subseteq \widehat{R_M}$, it follows that $Q \cap \widehat{R}$ is a nonzero prime ideal of $\widehat{R}$. Since $\widehat{R}$ is a Krull domain, there is some height-one prime ideal $P$ of $\widehat{R}$ such that $P \subseteq Q \cap \widehat{R}$. We conclude that $0 \neq P \cap R \subseteq Q \cap R = Q \cap R_M \cap R = M_M \cap R = M$, and so $P \cap R = M$. In particular, for each maximal ideal $N$ of $R$, there is a height-one prime ideal $N^{\prime}$ of $\widehat{R}$ such that $N^{\prime} \cap R = N$. Since each nonzero element of $\widehat{R}$ is contained in only finitely many height-one prime ideals of $\widehat{R}$, we have that $R$ is of finite character. It follows by Proposition \ref{prop22} that $R$ is an FSD.
\epf

\section{RINGS OF POLYNOMIALS} \label{S3}

We now turn our attention to investigating the FMD property as it pertains to certain rings of polynomials. Motivated by the original algebraic geometric context for primary decompositions, we find complete information with respect to the ring $k[X^2,X^3]$, that is, the ring of all polynomials over the field $k$ that lack a linear term, in Theorem \ref{theo10}. We then focus on standard polynomial rings $R[X]$, where $R$ is a domain and $X$ an indeterminate over $R$, and discover in Corollary \ref{coro14.1} that such a ring is an FMD when $R$ is a special type of Dedekind domain.

\bt \label{theo10}
Let $k$ be a field. Then the following are equivalent:

(1) $k[X^2,X^3]$ is an FMD;

(2) $k[X^2,X^3]$ is an FSD;

(3) $k$ is finite.
\et

\bpf
(1) $\Leftrightarrow$ (2): Since every FSD in an FMD by Corollary \ref{coro1.2}, we immediately have that (2) $\Rightarrow$ (1). Conversely, suppose (1). Clearly, $k[X^2,X^3]$ is a Noetherian one-dimensional domain. Consequently, $k[X^2,X^3]$ is an FSD by Theorem \ref{theo23}.

(2) $\Leftrightarrow$ (3): Suppose (2). Then the ideal $(X^4)$ is contained in only finitely many ideals of $k[X^2,X^3]$, whence the set $\{(X^2+bX^3,X^4) \mid b \in k\}$ is finite. It is sufficient then to show that there is a bijection between $k$ and $\{(X^2+bX^3,X^4) \mid b \in k\}$.

To see this, let $a,b \in k$ and suppose that $(X^2 + bX^3, X^4) = (X^2 + aX^3, X^4)$. Then for some $f, g \in k[X^2,X^3]$ it must be the case that $X^2 + bX^3 = f X^4 + g X^2 + g aX^3 \Rightarrow X^2(1 - g - f X^2) = (g a - b)X^3 \Rightarrow X^2 \,|\, (g a - b)X^3$. So, the $X^3$ term of $(g a - b)X^3$ must be $0$. Thus, $g_0 a = b$, where $g_0$ is the constant term of $g$. However, $X^3 \,|\, (1 - g)X^2 - f X^4$ and so it must be the case that the $X^2$ term is $0$. Thus, $g_0 = 1$, and so $a = b$. The bijection is thus established.

Conversely, suppose (3). Put $T = k[X]$ and put $R = k[X^2,X^3]$. Then $T$ is clearly integral over $R$ and $T$ is a finitely generated $R$-module. Since $T$ has the finite norm property, $R$ is an FSD by \cite[Theorem 3.6]{HL0}.
\epf

We next set out to investigate which ideals $I$ of a general polynomial domain $R[X]$ are molecules, not only towards a general understanding of such matters, but to discover sufficient conditions for a polynomial domain to be an FMD. Proposition \ref{prop11} and its associated Corollary \ref{coro11.1} below reveal that the context of UFD's allows for some definitive information along these lines.

Recall that for an ideal $I$ of the domain $R$, the {\em t}-closure of $I$ is given by $I_t = \cup \{J_v \mid J \subseteq I \textrm{ is a finitely generated ideal of } R\}$, where $J_v = (J^{-1})^{-1}$. An ideal is called a {\em t}-ideal if it coincides with its {\em t}-closure. Note that every nonzero principal ideal is a {\em t}-ideal and, moreover, the domain $R$ is a UFD if and only if every {\em t}-ideal of $R$ is principal. Furthermore, if $S \subseteq R\backslash\{0\}$ is a multiplicatively closed subset of $R$ and $J$ is a {\em t}-ideal of the quotient overring $S^{-1}R$, then $J \cap R$ is a {\em t}-ideal of $R$.

\bp \label{prop11}
Let $R$ be a UFD and $I$ a nonzero ideal of $R$. If $I \subsetneq I_t \subsetneq R$, then $I$ is a compound ideal of $R$.
\ep

\bpf
Let $I \subsetneq I_t \subsetneq R$. Then $I$ is properly contained in the proper multiplication ideal $I_t$, and thus $I$ is not a molecule of $R$ by Proposition \ref{prop0}.
\epf

\bc \label{coro11.1}
Let $R$ be a UFD and $I$ a nonprincipal ideal of $R[X]$. If $I \cap R = 0$, then $I$ is a compound ideal.
\ec

\bpf
Let $I \cap R = 0$ and let $K$ be the quotient field of $R$. Observe that $R[X]$ is a UFD and $IK[X] \cap R[X] \subsetneq R[X]$. Since $IK[X]$ is a principal ideal of $K[X]$ and $K[X]$ is a quotient overring of $R[X]$, we have that $IK[X] \cap R[X]$ is a {\em t}-ideal of $R[X]$. This implies that $I_t \subseteq (IK[X] \cap R[X])_t = IK[X] \cap R[X]$. We conclude that $I \subsetneq I_t \subsetneq R[X]$, and hence $I$ is not a molecule by Proposition \ref{prop11}.
\epf

While Proposition \ref{prop11} above established that certain ideals in a UFD are compound, Proposition \ref{prop12} below provides for a wealth of ideals in an arbitrary domain that are molecules.

\bp \label{prop12}
Let $R$ be a domain and let $P$ and $Q$ be nonzero locally principal ideals of $R$. If $P$ is a prime ideal of $R$ and $P+Q$ is a maximal ideal of $R$, then $P + Q^n$ is a molecule of $R$ for each $n \in \mathbb{N}$.
\ep

\bpf
Let $P$ be a prime ideal of $R$, $P+Q$ a maximal ideal of $R$, and $n \in \mathbb{N}$. Put $M = P+Q$. If $Q \subseteq P$, then $M = P +Q^n = P$, and thus $P+Q^n$ is a molecule of $R$ by Proposition \ref{prop6.1}(c). Therefore, we can assume that $Q \nsubseteq P$. Suppose to the contrary that $P+Q^n$ is not a molecule of $R$. Then there are proper ideals $J$ and $K$ of $R$ for which $P + Q^n = JK$. It follows that $M = \sqrt{P+Q^n} = \sqrt{J} \cap \sqrt{K}$, and hence $\sqrt{J} = \sqrt{K} = M$. As such, $J$ and $K$ are contained in $M$. We conclude that $P \subseteq P +Q^n = JK \subseteq M^2 = P^2 +PQ + Q^2$. However, since $P$ and $Q$ are both locally principal, there exist $p,h \in R_M$ such that $P_M = pR_M$ and $Q_M = hR_M$. Thus, we have that $pR_M = P_M \subseteq P_M^2 + P_MQ_M + Q_M^2 = p^2R_M + phR_M + h^2R_M$, hence there are $g_1,g_2,g_3 \in R_M$ such that $p = p^2g_1 + phg_2 + h^2g_3$. Observe that $h^2g_3 \in P_M$. If $h \in P_M$, then $Q \subseteq Q_M \cap R \subseteq P_M \cap R = P$, a contradiction. We may conclude that $h \not\in P_M$, and so $g_3 = pg_4$ for some $g_4 \in R_M$. But this implies that $1 = pg_1 + hg_2 + h^2g_4 \in P_M + Q_M = M_M$, a contradiction. Therefore, $P +Q^n$ is a molecule of $R$.
\epf

Theorem \ref{theo13} reveals that if $R$ is further assumed to be a Dedekind domain, a characterization along the lines of Proposition \ref{prop12} is available.

\bt \label{theo13}
Let $R$ be a Dedekind domain, $P$ a nonzero prime ideal of $R$, and $I$ a proper ideal of $R[X]$ such that $I \supsetneq P[X]$. Then $I$ is a molecule of $R[X]$ if and only if there are some $n \in \mathbb{N}$ and $f \in R[X]$ such that $P[X] + fR[X]$ is a maximal ideal of $R[X]$ and $I = P[X] + f^nR[X]$.
\et

\bpf
Observe that $R[X]$ is a two-dimensional Noetherian domain and $P[X]$ is an invertible height-one prime ideal of $R[X]$. Therefore, every prime ideal of $R[X]$ that contains $I$ is a maximal ideal of $R[X]$.

Suppose that $I$ is a molecule of $R[X]$. It follows from Corollary \ref{coro7.12} that $I$ is primary. Put $M = \sqrt{I}$. Then $M$ is a maximal ideal of $R[X]$. Note that $R[X]/P[X] \cong (R/P)[X]$ is a PID and $M/P[X]$ is a maximal ideal of $R[X]/P[X]$. Consequently, $M = P[X] + fR[X]$ for some $f \in M$. Moreover, $\sqrt{I/P[X]} = M/P[X]$, and so $I/P[X] = (M/P[X])^n$ for some $n \in \mathbb{N}$. Therefore, $I = P[X] + M^n = P[X] + f^nR[X]$, as desired.

Since the converse is an immediate consequence of Proposition \ref{prop12}, the proof is complete.
\epf

We pause briefly to note through Proposition \ref{prop13.3} that, even if $R$ is a PID, there exist nonprincipal molecules of $R[X]$ other than those of the form $P[X]+f^n[X]$, where $P$ is a nonzero prime ideal of $R$ and $P[X] + fR[X]$ is a maximal ideal of $R[X]$. However, as the proof of Proposition \ref{prop13.3} makes clear, demonstrating that such ideals are molecules can prove to be relatively involved.

\bp \label{prop13.3}
Let $R$ be a PID and $p$ a prime element of $R$. The ideal $(X^2,p^2)$ in $R[X]$ is a molecule.
\ep

\bpf
Suppose to the contrary that $P = (X^2,p^2)$ is compound, so that $P = IJ$ for proper ideals $I$ and $J$. Then $P = (P:J)J$, where $P:J$ is also proper, owing to the fact that $R[X]$ is a Noetherian domain. Let $M$ be the maximal ideal $(X,p)$ in $R[X]$ and observe that $M^3 \subseteq P$, so that $P:J$ and $J$ are both contained in $M$. We will show that if $P:J \neq M$, then $J = M$.

Suppose that $P:J \neq M$. Since $M$ is maximal, it must be the case that $MJ \nsubseteq P$. Since $MJ \subseteq M^2 = (X^2, pX, p^2)$ and $MJ \nsubseteq P = (X^2, p^2)$, there exists $n \in R$ such that $\gcd(n,p) = 1$ and $npX \in MJ-P$. Let $u, v \in R$ be such that $1 = un + vp$. Then $pX = u(npX) + vp^2X \in MJ$. Thus, $MJ = M^2$.

Next, we show that $J = M$. Note that it suffices to show that $p, X \in J$. Since $pX \in MJ = XJ + pJ$, there exist $\alpha, \beta \in J$ such that $pX = X\alpha + p\beta$. Since $p \,|\, \alpha$ and $X \,|\, \beta$, it follows that $\alpha = hp$ and $\beta = gX$ for some $h,g \in R[X]$. Thus, $pX = hpX + gXp = pX(h + g)$, whence $h + g = 1$. However, since $h = b + Xf$ for some $b \in R$ and $f \in R[X]$, it must be the case that $bp + pXf = \alpha \in J$, and so $bp \in J$. Moreover, $(1 - b)X - X^2f = \beta \in J$, and hence $(1 - b)X \in J$. We then have the following two cases:

\textit{Case I}. Suppose that $p \mid_R b$. Then $bX \in J$. Thus, $X = bX + (1 - b)X \in J$. Since $p^2 \in XJ + pJ$, there exist $y,z \in J$ such that $p^2 = Xy + pz$. Observe that $p \,|\, y$, so that $y = pw$ for some $w \in R[X]$. It follows that $p = Xw + z \in J$.

\textit{Case II}. Suppose that $p \nmid_R b$. Then $pR + bR = R$, and thus $p \in pR[X] = p^2R[X] + bpR[X] \subseteq J$. Since $X^2 \in XJ + pJ$, there exist $y,z \in J$ such that $X^2 = Xy + pz$. Note that $X \,|\, z$, and thus $z = Xw$ for some $w \in R[X]$. This implies that $X = y + pw \in J$.

Therefore, either $P:J = M$ or $J = M$. It then follows that $M$ divides $P$, and so $P = (P:M)M$. Since $M^3 \subseteq P$, we have that $M^2 \subseteq P:M$. Let $f \in P:M$. Then there exist $\alpha, \beta, \gamma, \delta \in R[X]$ for which $fX = \alpha X^2 + \beta p^2$ and $fp = \gamma X^2 + \delta p^2$. So, $f = \alpha X + \frac{\beta}{X} p^2$ and $f = \frac{\gamma}{p}X^2 + \delta p$. Thus, $f \in (X,p^2) \cap (X^2,p) = M^2$. Since $f$ was an arbitrary element of $P:M$, it follows that $P:M = M^2$. However, this implies that $P = (P:M)M = M^3$, a contradiction. Therefore, $P$ is a molecule.
\epf

We now come to the point where we can present a subclass of Dedekind domains $R$ for which $R[X]$ is an FMD. This task is formalized in Corollary \ref{coro14.1} below as a consequence to a theorem (Theorem \ref{theo14}) revealing how information related to the prime ideals of a domain is enough to guarantee that the domain is an FMD.

We first provide a lemma that highlights a particular behavior of invertible ideals that can be utilized in the context of Noetherian domains.

\bl \label{lemm13.5}
Let $R$ be a domain that satisfies the ascending chain condition on invertible ideals and let $I$ be a nonzero ideal of $R$. Then $I$ is contained in an invertible ideal of $R$ that is minimal amongst all invertible ideals containing $I$.
\el

\bpf
Suppose that $I$ is a nonzero ideal of the domain $R$ which is not contained in an invertible ideal of $R$ that is minimal amongst all invertible ideals containing $I$. Then there exists a properly descending chain $R \supsetneq J_1 \supsetneq J_2 \supsetneq \, \cdots$ of invertible ideals that contain $I$. Choose some $0 \neq a \in I$ and put $K_i = aJ_i^{-1}$ for $i =1,2, \, \ldots$, so that each $K_i$ is necessarily an invertible ideal of $R$. Moreover, $K_1 \subsetneq K_2 \subsetneq K_3 \subsetneq \, \cdots$. Therefore, $R$ does not satisfy the ascending chain condition on invertible ideals, and the result follows.
\epf

\bt \label{theo14}
Let $R$ be a two-dimensional Noetherian domain such that every height-one prime ideal of $R$ is invertible and $R/M$ is finite for every height-two prime ideal $M$ of $R$. Then $R$ is an FMD.
\et

\bpf
First we show that if $I$ is an ideal of $R$ that is not contained in a height-one prime ideal of $R$, then $I$ is contained in only finitely many ideals of $R$.

Let $I$ be an ideal of $R$ that is not contained in a height-one prime ideal of $R$. Observe that if $M$ is a height-two prime ideal of $R$ and $n \in \mathbb{N}$, then $R/M^n$ is finite (since $M^k/M^{k+1}$ is a finite-dimensional vector space over $R/M$ for each $k \in \mathbb{W}$ and $\left|R/M^n\right| = \prod_{k=0}^{n - 1}\left|M^k/M^{k+1}\right| < \infty$). Let $V(I)$ be the set of prime ideals of $R$ containing $I$. Then $V(I)$ consists only of maximal ideals of $R$. Since $R$ is Noetherian, $V(I)$ is finite, and furthermore the elements of $V(I)$ are pairwise comaximal. Hence $\sqrt{I} = \prod_{M \in V(I)} M$. Again owing to the fact that $R$ is Noetherian, there is some $n \in \mathbb{N}$ such that $J := \prod_{M \in V(I)} M^n \subseteq I$. It follows by the Chinese remainder theorem that $R/J \cong \prod_{M \in V(I)} R/M^n$. Consequently, $R/J$ is finite, and thus $R/I$ is finite. Therefore, $I$ is contained in only finitely many ideals of $R$.

By Corollary \ref{coro1.1}, it is sufficient to show that every nonzero ideal of $R$ is divisible by only finitely many molecules of $R$. Let $I$ be a nonzero ideal of $R$. By Lemma \ref{lemm13.5}, there is an invertible ideal $J$ of $R$ that contains $I$ and is minimal amongst all of the invertible ideals of $R$ that contain $I$. Assume that $IJ^{-1}$ is contained in a height-one prime ideal $P$ of $R$. By assumption, $P$ must be invertible. Then $I \subseteq PJ \subseteq J$, and hence $PJ = J$, by minimality of $J$. However, this implies that $P = R$, a contradiction. Therefore, $IJ^{-1}$ is contained in only finitely many ideals of $R$ by the above work.

Now, let $K$ be a molecule of $R$ that divides $I$. Then $I = KA$ for some ideal $A$ of $R$. Since $R$ is Noetherian, the set of height-one prime ideals of $R$ that contain $I$ is finite. Thus, it is sufficient to show that $K$ is either a height-one prime ideal of $R$ or $IJ^{-1} \subseteq K$.

Suppose that $K$ is contained in a height-one prime ideal $P$ of $R$. Since $P$ is invertible by hypothesis, $K = PB$ for some ideal $B$ of $R$, and hence $K = P$.

On other other hand, suppose $K$ is not contained in a height-one prime ideal of $R$. Note that there is an invertible ideal $B$ of $R$ that contains $J + A$ and is minimal amongst all the invertible ideals of $R$ that contain $J + A$. Assume that $J \neq B$. Then $JB^{-1}$ is a proper invertible ideal of $R$ and thus it is contained in a height-one prime ideal $P$ of $R$. But then $KAB^{-1} = IB^{-1} \subseteq JB^{-1} \subseteq P$. We conclude that $AB^{-1} \subseteq P$, and hence $J + A \subseteq BP \subseteq B$. But, again by minimality of $B$, it follows that $BP = B$, and so $P = R$, a contradiction. Therefore, $A \subseteq B = J$, hence $IJ^{-1} = KAJ^{-1} \subseteq K$, as desired.
\epf

As in \cite{CL} and \cite{LM}, a domain $R$ is said to have the \textit{finite norm property} if $R/I$ is finite for every nonzero ideal $I$ of $R$.

\bc \label{coro14.1}
(a) If $R$ is a domain such that $R[X]$ is an FMD, then $R$ is an FMD.

(b) If $R$ is a Dedekind domain with the finite norm property, then $R[X]$ is an FMD.
\ec

\bpf
(a) Since $R \subseteq R[X]$ is a survival extension, it is an immediate consequence of Proposition \ref{prop4.51} and Corollary \ref{coro1.01} that $R$ has unit-cancellation for ideals. Let $I$ be a nonzero ideal of $R$. Let $\mathcal{I}$ be the set of ideals of $R$ that divide $I$ and $\mathcal{J}$ be the set of ideals of $R[X]$ that divide $I[X]$. Observe that $f$:$\,\, \mathcal{I} \to \mathcal{J}$ given by $f(L) = L[X]$ is a well-defined injective map. It follows by Corollary \ref{coro1.01} that $\mathcal{J}$ is finite, and thus $\mathcal{I}$ is finite. We conclude from Corollary \ref{coro1.01} that $R$ is an FMD.

(b) Let $K$ be the quotient field of $R$. If $R=K$, then $R[X]$ is a PID, and so $R[X]$ is an FMD. Thus, we may assume that $R \neq K$. As such, it is well-known that $R[X]$ is a two-dimensional Noetherian domain.

We first show that every height-one prime ideal of $R[X]$ is invertible. Let $P$ be a height-one prime ideal of $R[X]$. We then have the following two cases:

\textit{Case I.} Suppose that $P \cap R \neq 0$. Put $Q = P \cap R$. Then $Q$ is a nonzero prime ideal of $R$ and $P = Q[X]$. Since $Q$ is an invertible ideal of $R$ and $P$ is an extension ideal of $Q$, we have that $P$ is an invertible ideal of $R[X]$.

\textit{Case II.} Suppose that $P \cap R = 0$. Then there is some nonzero $f \in R[X]$ for which $P = fK[X] \cap R[X]$. Let $c(f)$ be the ideal of $R$ generated by the coefficients of $f$ (that is, $c(f)$ is the content of $f$). Since $R$ is integrally closed, it follows by \cite[Lemme 1]{Q} that $P = fc(f)^{-1}[X]$. Since $c(f)$ is an invertible ideal of $R$, we have that $c(f)^{-1}$ is an invertible fractional ideal of $R$. Therefore, $c(f)^{-1}[X]$ is an invertible fractional ideal of $R[X]$, and so $P$ is an invertible ideal of $R[X]$.

Now, we show that $R[X]/M$ is finite for all height-two prime ideals of $R[X]$. Let $M$ be a height-two prime ideal of $R[X]$. Put $Q = M \cap R$. Since every upper to zero has height-one, it follows that $Q$ is a nonzero prime ideal of $R$. In addition, since $Q[X] \subsetneq Q + XR[X]$ and $Q[X]$ and $Q + XR[X]$ are prime ideals of $R[X]$, it follows that $Q[X]$ is a height-one prime ideal of $R[X]$. Moreover, there is a ring isomorphism $\phi$:$\,\,R[X]/Q[X] \to (R/Q)[X]$. Put $N = \phi(M/Q[X])$. Since $M$ is a maximal ideal of $R[X]$, it follows that $N$ is a maximal ideal of $(R/Q)[X]$. However, since $R$ has the finite norm property, $R/Q$ is a finite field, whence $(R/Q)[X]$ has the finite norm property. Therefore, $R[X]/M \cong (R/Q)[X]/N$ is finite, as desired.

We conclude that $R[X]$ is an FMD by an application of Theorem \ref{theo14}.
\epf

\bc \label{coro14.2}
If $\mathcal{O}_K$ is the ring of integers of the algebraic number field $K$, then $\mathcal{O}_K[X]$ is an FMD.
\ec

\bpf
This follows immediately from Corollary \ref{coro14.1} and the well-known fact that the ring of integers in any algebraic number field is a Dedekind domain with the finite norm property.
\epf

\section{FMD'S AND THE D+M CONSTRUCTION} \label{S4}

We conclude this paper with a brief section on the intersection of the FMD property and the classical $D+M$ construction. An exploration of the ``$*$-factorability" of $D+M$ construction was conducted in \cite{AKP}, where $*$ is a star operation on the fractional ideals of the domain. For our purposes here, we limit the discussion to $* = d$ and, inspired by \cite[Proposition 7, Corollary 8]{AKP}, first provide a characterization for when $D+M$ is a molecular domain.

\bt \label{theo15}
Let $V$ be a valuation domain of the form $V = K + M$, where $M$ is the nonzero maximal ideal of $V$ and $K$ is a field. Let $D$ be a proper subring of $K$, and put $R = D + M$. Then $R$ is a molecular domain if and only if $V$ is a DVR and $D$ is a field.
\et

\bpf
First, suppose that $R$ is molecular. Then $D$ is a field by \cite[Proposition 7]{AKP}. Now, suppose that $V$ is not a DVR. Then there exists a prime ideal $P$ of $V$ that is not finitely generated. Thus, $PM = P$. However, then $P$ is a prime ideal of $R$ for which $P = PM$ in $R$, as well. This implies that $P$ is a compound prime ideal in $R$, contradicting the fact that prime ideals are molecules in a molecular domain. Therefore, $V$ must be a DVR.

Conversely, suppose that $V$ is a DVR and $D$ is a field. Let $a \in V$ be such that $M = aV$. We claim that every nonzero proper ideal of $R$ has the form $a^nF + a^{n+1}V$ for some nonzero $D$-subspace $F$ of $K$ and $n \geq 1$ (cf. \cite[Theorem 2.1(k)]{BG}). Let $I$ be a nonzero proper ideal of $R$. Since $V$ is a DVR, there exists $b \in V$ for which $IV = bV$ and furthermore it must be the case that $I = Fb + Mb$ in $R$ for some $D$-subspace $F$ of $K$ \cite[Theorem 2.1(n)]{BG}. Since every nonzero proper ideal of $V$ is a power of $M$, one may take $b = a^n$ for some $n \geq 1$. Hence, $I = Fa^n + a^{n+1}V$. Moreover, if $F = 0$, then $I = a^{n+1}V = a^{n+1}K + a^{n+2}V$, where $K$ is necessarily nonzero. Thus, $F$ can be assumed to be nonzero, and the claim is established.

Now, it is straightforward to verify that for a nonzero proper ideal $I$ of $R$, it is the case that $I = a^nF + a^{n+1}V = (aR)^{n-1}(aF + a^2V)$. One can further show that $M^2 \subsetneq aF + a^2V$. Thus, $aF + a^2V$ is a molecule by \cite[Proposition 2.2]{HL1}. Observe that since $D$ is itself a $D$-subspace of $K$, then $aR$ must also be a molecule (cf. \cite[Example 4]{AKP}). Thus, $(aR)^{n-1}(aF + a^2V)$ is a molecularization for $I$. Therefore, $R$ is a molecular domain.
\epf

Our final result, Theorem \ref{theo16}, shows that the addition of the simple condition ``$K$ is finite" (where $K$ is the residue field of the valuation domain $V$) to the conditions in Theorem \ref{theo15} above is all that is needed to characterize when $D+M$ is an FMD.

\bt \label{theo16}
Let $V$ be a valuation domain of the form $K + M$, where $M$ is the nonzero maximal ideal of $V$, and $K$ is a field. Let $D$ be a proper subring of $K$, and put $R = D + M$. The following are equivalent:

(1) $R$ is an FSD;

(2) $V$ is a DVR, $D$ is a field, and $K$ is finite;

(3) $V$ is a DVR, $D$ is a field, and $K^*/D^*$ is finite;

(4) $R$ is an FFD;

(5) $R$ is an FMD.
\et

\bpf
The equivalence of conditions (1), (2), (3), and (4) was established in \cite[Theorem 4.3]{HL0}. However, Corollary \ref{coro1.2} shows that in general (1) $\Rightarrow$ (5) and Corollary \ref{coro1.3} shows that in general (5) $\Rightarrow$ (4). The proof is thus complete.
\epf

\end{document}